# Exact solution of Helmholtz equation
# for the case of non-paraxial Gaussian beams.

**Sergey V. Ershkov**

Institute for Time Nature Explorations,

M.V. Lomonosov's Moscow State University,

Leninskie gory, 1-12, Moscow 119991, Russia

e-mail: sergej-ershkov@yandex.ru

A new type of exact solutions of the full 3 dimensional *spatial* Helmholtz equation for the case of non-paraxial Gaussian beams is presented here.

We consider appropriate representation of the solution for Gaussian beams *in a spherical coordinate system* by substituting it to the full 3 dimensional spatial Helmholtz Equation.

Analyzing the structure of the final equation, we obtain that governing equations for the components of our solution are represented by the proper *Riccati* equations of complex value, which has no analytical solution in general case.

But we find one of the possible exact solutions which is proved to satisfy to such an equations for Gaussian beams.

## 1. Introduction.

The full 3-dimensional *spatial* Helmholtz equation provides solutions that describe the propagation of waves over space (*e.g., electromagnetic waves*) under a proper boundary conditions; it should be presented in a spherical coordinate system $R, \theta, \varphi$ as below [1-2]:

$$\Delta A + k^2 A = 0, \qquad (1.1)$$

- where $\Delta$ - is the Laplacian, $k$ is the wavenumber, and $A$ is the amplitude. So, the derivation advanced in this manuscript starts with the scalar Helmholtz equation expressed in spherical co-ordinates.

Besides, in spherical coordinate system [3]:

$$\Delta A = \frac{\partial^2 A}{\partial R^2} + \frac{2}{R}\frac{\partial A}{\partial R} + \frac{1}{R^2 \sin^2\theta}\frac{\partial^2 A}{\partial \varphi^2} + \frac{1}{R^2}\frac{\partial^2 A}{\partial \theta^2} + \frac{1}{R^2}\cot\theta \frac{\partial A}{\partial \theta} \ .$$

Special solutions to this equation have generated continuing interest in the optical physics community since the discovery of unusual non-diffracting waves such as Bessel and Airy beams [4-6].

Let us search for solutions of Eq. (1.1) in a *classical* form of Gaussian beams [7-9], which could be presented in Cartesian coordinate system $X, Y, Z$ as below [10]:

$$A = a \cdot \frac{w_0}{w(Z)} \exp\left[-\frac{X^2+Y^2}{w^2(Z)} - ikZ - ik\frac{X^2+Y^2}{2r(Z)} + i\zeta(Z)\right]$$

- where $w(Z), r(Z), \zeta(Z)$ – are the real functions, describing the appropriate parameters of a beam; $w(Z)$ is the beam waist size, $r(Z)$ is a wavefront radius of curvature and $\zeta(Z)$ is the *Gouy's phase shift* properly [10].

The *classical* form of Gaussian beams above could be also represented as below

$$\exp\left[i\left(\zeta(Z) - kZ + i\cdot\ln w(Z) + \left(\frac{i}{w^2(Z)} - \frac{k}{2R(Z)}\right)\cdot(X^2+Y^2)\right)\right] = \exp\left[i\left(p(Z) + \frac{X^2+Y^2}{2q(Z)}\right)\right]$$

- here $p(Z)$ is the complex phase-shift of the waves during their propagation along the $Z$ axis; $q(Z)$ is the proper complex parameter of a beam, which is determining the Gaussian profile of a wave in the transverse plane at position $Z$.

Besides, let us also note that at the *left* part of the expression above we express the term $(1/w(Z))$ in a form for Gaussian beams, as $\exp(i^2\cdot\ln w(Z)) = \exp(-\ln w(Z))$. The *right* part of the expression above could be transformed in a *spherical* coordinate system to the form below:

$$A = a\cdot\exp\left[i\left(p(R,\theta) + \frac{R^2\cdot\sin^2\theta}{2q(R,\theta)}\right)\right] \qquad (*)$$

The solution (*) is additionally assumed to be independent of the azimuthal co-ordinate to observe it under well-known paraxial approximation [10] also.

Then having substituted the expression (*) into Eq. (1.1), we should obtain ($\theta \neq 0$):

$$\frac{\partial^2 p(R,\theta)}{\partial R^2} + \frac{\partial^2\left(\frac{R^2}{q(R,\theta)}\right)}{\partial R^2}\frac{\sin^2\theta}{2} + i\cdot\left(\frac{\partial p(R,\theta)}{\partial R} + \frac{\partial\left(\frac{R^2}{q(R,\theta)}\right)}{\partial R}\frac{\sin^2\theta}{2}\right)^2 + \frac{2}{R}\left(\frac{\partial p(R,\theta)}{\partial R} + \frac{\partial\left(\frac{R^2}{q(R,\theta)}\right)}{\partial R}\frac{\sin^2\theta}{2}\right) +$$

(1.2)

$$+\frac{1}{R^2}\cdot\frac{\partial^2 p(R,\theta)}{\partial\theta^2} + \frac{1}{2}\frac{\partial^2\left(\frac{\sin^2\theta}{q(R,\theta)}\right)}{\partial\theta^2} + \frac{i}{R^2}\cdot\left(\frac{\partial p(R,\theta)}{\partial\theta} + \frac{\partial\left(\frac{\sin^2\theta}{q(R,\theta)}\right)}{\partial\theta}\frac{R^2}{2}\right)^2 + \frac{\cot\theta}{R^2}\cdot\left(\frac{\partial p(R,\theta)}{\partial\theta} + \frac{\partial\left(\frac{\sin^2\theta}{q(R,\theta)}\right)}{\partial\theta}\frac{R^2}{2}\right) =$$

$$= i\cdot k^2.$$

## 2. Exact solutions.

Let us re-designate appropriate term in (*) as below:

$$f(R,\theta) = p(R,\theta) + \frac{R^2 \cdot \sin^2\theta}{2q(R,\theta)}.$$

In such a case, Eq. (1.2) could be transformed as below ($\theta \neq 0$):

$$\frac{\partial^2 f(R,\theta)}{\partial R^2} + i \cdot \left(\frac{\partial f(R,\theta)}{\partial R}\right)^2 + \frac{2}{R} \cdot \left(\frac{\partial f(R,\theta)}{\partial R}\right) + $$

$$+ \frac{1}{R^2} \cdot \left(\frac{\partial^2 f(R,\theta)}{\partial \theta^2} + i \cdot \left(\frac{\partial f(R,\theta)}{\partial \theta}\right)^2 + \cot\theta \cdot \left(\frac{\partial f(R,\theta)}{\partial \theta}\right)\right) - i \cdot k^2 = 0$$

(2.1)

Thus, all possible solutions for representing of Gaussian beams in a form (*) are described by the Equation (2.1).

But we should especially note that during the process of the obtaining of a solution (for example, if we are simply assuming a special eikonal solution [10-11] to the Helmholtz equation), some of main features of the solution could be reduced; so, such a solution need not have any relation to a Gaussian form (*).

Besides, one of the obvious solutions of PDE-equations (2.1):

$$f(R, \theta) = f_1(R) + f_2(\theta) \qquad (**)$$

- where $f_1(R)$, $f_2(\theta)$ – are the functions of *complex* value.

Let us assume as below:

$$\frac{\partial^2 f(R,\theta)}{\partial \theta^2} + i\cdot\left(\frac{\partial f(R,\theta)}{\partial \theta}\right)^2 + \cot\theta\cdot\left(\frac{\partial f(R,\theta)}{\partial \theta}\right) = C \qquad (2.2)$$

- here C – is a constant of *complex* value. For such a case, Eq. (2.1) could be reduced as below ($\theta \neq 0$):

$$\frac{\partial^2 f(R,\theta)}{\partial R^2} + i\cdot\left(\frac{\partial f(R,\theta)}{\partial R}\right)^2 + \frac{2}{R}\cdot\left(\frac{\partial f(R,\theta)}{\partial R}\right) + \frac{C}{R^2} - i\cdot k^2 = 0 \qquad (2.3)$$

## 3. Presentation of exact solution.

Under assumption (**), Eq. (2.2) could be represented as below:

$$\left(\frac{d f_2}{d\theta}\right) = y(\theta) \quad \Rightarrow \quad y'(\theta) = -i\cdot y^2 - \cot\theta \cdot y + C,$$

(3.1)

$$y(\theta) = \csc\theta \cdot u(\theta) \quad \Rightarrow \quad u'(\theta) = -(i\cdot\csc\theta)\cdot u^2 + C\cdot\sin\theta,$$

- where the last equation is known to be the *Riccati* ODE [3], which has no solution in general case. But if $C = 0$, Eq. (3.1) has a proper solution ($C_0 = $ const):

$$u'(\theta) = -(i\cdot\csc\theta)\cdot u^2, \qquad u(\theta) = \frac{1}{\left(C_0 + i\cdot\int\csc\theta\, d\theta\right)} \Rightarrow$$

(3.2)

$$\frac{d f_2}{d\theta} = \frac{\csc\theta}{\left(C_0 + i\cdot\int\csc\theta\, d\theta\right)} \quad (C_0 = 0) \Rightarrow f_2 = -i\cdot\ln\left(\int\csc\theta\, d\theta\right)$$

Besides, Eq. (2.3) could be presented as below ($C = 0$):

$$\left(\frac{df_1}{dR}\right) = y_1(R) \quad \Rightarrow \quad y_1'(R) = -i \cdot y_1^2 - \frac{2}{R}y_1 - \left(\frac{C}{R^2} - i \cdot k^2\right), \tag{3.3}$$

$$f_1(R) = \int y_1(R)\, dR\ .$$

- where the last *Riccati* ODE (3.3) has a proper solution below if $C = 0$ (see [3], the case 1.104).

Indeed, let us assume ($k \neq 0$, $R \neq 0$):

$$y_1 = u_1 + \frac{i}{R}, \quad y_1'(R) = -i \cdot y_1^2 - \frac{2}{R}y_1 + i \cdot k^2 \Rightarrow$$

$$\Rightarrow u_1'(R) = -i \cdot u_1^2 + i \cdot k^2 \Rightarrow \int \frac{du_1}{k^2 - u_1^2} = i \cdot R$$

$$\Rightarrow \begin{cases} u_1 = k \cdot \tanh(i \cdot k \cdot R),\ |i \cdot \tan(k \cdot R)| < 1, \\ \\ u_1 = k \cdot \coth(i \cdot k \cdot R),\ |i \cdot \tan(k \cdot R)| > 1, \end{cases}$$

- then, we obtain:

$$\begin{cases} f_1 = -i \cdot \ln\cosh(i \cdot k \cdot R) + i \cdot \ln R,\ |k \cdot R| < \pi/4, \\ \\ f_1 = -i \cdot \ln\sinh(i \cdot k \cdot R) + i \cdot \ln R,\ |k \cdot R| > \pi/4. \end{cases} \tag{3.4}$$

Taking into consideration the expression (**) for the solution as well as (3.2)-(3.4), let us finally present a new type of *non-paraxial* solution, which is proved to satisfy to the Helmholtz equation (1.1), as below:

$$\begin{cases} A = a \cdot \left( \int \csc\theta \, d\theta \right) \cdot \dfrac{\cosh(i \cdot k \cdot R)}{R}, & |k \cdot R| < \pi/4, \\[2ex] A = a \cdot \left( \int \csc\theta \, d\theta \right) \cdot \dfrac{\sinh(i \cdot k \cdot R)}{R}, & |k \cdot R| > \pi/4, \end{cases}$$

- or

$$\begin{cases} A = (k \cdot a) \cdot \ln\left(\tan\left(\dfrac{\theta}{2}\right)\right) \cdot \dfrac{\cos(k \cdot R)}{k \cdot R}, & |k \cdot R| < \pi/4, \\[2ex] A = (k \cdot a) \cdot \ln\left(\tan\left(\dfrac{\theta}{2}\right)\right) \cdot \dfrac{i \cdot \sin(k \cdot R)}{k \cdot R}, & |k \cdot R| > \pi/4, \end{cases} \quad (3.5)$$

- where $\theta \in (0, \pi)$.

## **4. Discussions & conclusion.**

A new type of exact solutions of the full 3 dimensional *spatial* Helmholtz equation for the case of non-paraxial Gaussian beams is presented here.

We consider appropriate representation of the solution for Gaussian beams *in a spherical coordinate system* by substituting it to the full 3 dimensional spatial Helmholtz Equation.

Analyzing the structure of the final equation, we obtain that governing equations for the components of our solution are represented by the proper *Riccati* equations of complex value, which has no analytical solution in general case. We should note that a modern method exists for obtaining of the numerical solution of *Riccati* equations with a good approximation [12].

But we find one of the possible exact solutions (3.5) which is proved to satisfy to the Helmholtz equation (1.1) for beams (*).

Indeed, since the functions $g(R) = (\sin(k \cdot R))/(k \cdot R)$ or $g(R) = (\cos(k \cdot R))/(k \cdot R)$ in (3.5) are itself an exact solutions of the full Helmholtz equation (1.1) [1-2], the formula for the Laplacian in spherical coordinates gives for $A = h(\theta) \cdot g(R)$, $h(\theta) = \ln(tg(\theta/2))$:

$$\frac{1}{R^2}\frac{\partial^2 A}{\partial \theta^2} + \frac{1}{R^2}\cot\theta\frac{\partial A}{\partial \theta} = 0 \quad \rightarrow \quad \frac{d\left(\frac{1}{\sin\theta}\right)}{d\theta} + \cot\theta \cdot \left(\frac{1}{\sin\theta}\right) = 0,$$

- which is obviously valid for the range of parameter $\theta \in (0, \pi)$.

As for the appropriate example of *paraxial* approximation for such a *non-paraxial* exact solution (3.5) of the full Helmholtz equation (1.1), it could be easily obtained in the case $\theta \rightarrow +0$ (see the expression (3.5) above).

Let us express the real part of solution (3.5) in the Cartesian co-ordinates $X$, $Y$, $Z$ as below ($|k \cdot R| < \pi/4$, $r(X, Y) \ll Z$):

$$R = \sqrt{X^2 + Y^2 + Z^2} = \sqrt{(r(X,Y))^2 + Z^2}, \quad \cos\theta = \frac{Z}{\sqrt{(r(X,Y))^2 + Z^2}},$$

$$\tan\left(\frac{\theta}{2}\right) = \frac{1-\cos\theta}{\sqrt{1-(\cos\theta)^2}} = \sqrt{\frac{1-\cos\theta}{1+\cos\theta}} = \sqrt{\frac{\sqrt{(r(X,Y))^2 + Z^2} - Z}{\sqrt{(r(X,Y))^2 + Z^2} + Z}},$$

$$A = (k \cdot a) \cdot \ln\left(\tan\left(\frac{\theta}{2}\right)\right) \cdot \frac{\cos(k \cdot R)}{k \cdot R} = \quad (4.1)$$

$$= \frac{(k \cdot a)}{2} \cdot \ln\left(\frac{\sqrt{(r(X,Y))^2 + Z^2} - Z}{\sqrt{(r(X,Y))^2 + Z^2} + Z}\right) \cdot \frac{\cos(k \cdot \sqrt{(r(X,Y))^2 + Z^2})}{k \cdot \sqrt{(r(X,Y))^2 + Z^2}} \cong (k \cdot a) \cdot \ln\left(\frac{r(X,Y)}{2Z}\right) \cdot \frac{\cos(k \cdot Z)}{k \cdot Z}$$

As we know, a spherical-wave solution $g(R) = (\cos(k \cdot R))/(k \cdot R)$ could be *schematically* imagined in the Cartesian co-ordinate system as below [1-2]:

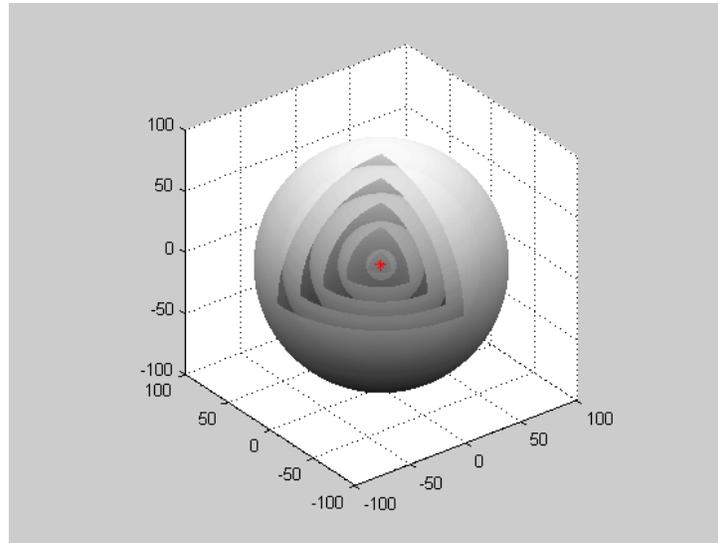

Fig. 1. A *schematic* plot of a spherical-wave type of the solutions.

- where each of spherical waves is assumed to be a concentric sphere evenly enlarging from a fixed point (a source of waves), see Fig.1-2.

The solution (3.5) differs from the spherical-waves on a factor $\ln(\tan(\theta/2))$, but the total energy of a beam should not exceed the total energy of the appropriate spherical-waves solution of Helmholtz equation. The energy of the beam is, of course, essentially the absolute magnitude of the solution spherically integrated over space.

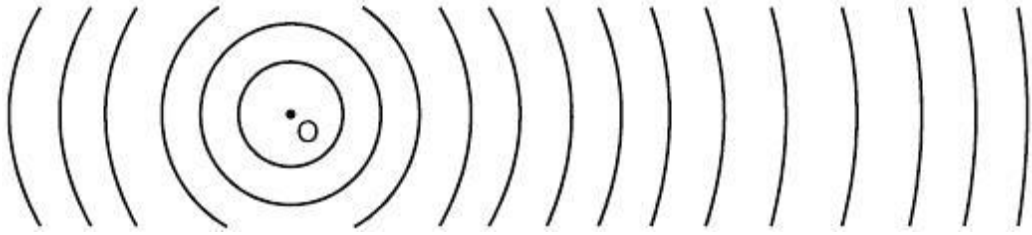

Fig. 2. A *schematic* plot of the plane spherical-wave solutions.

So, we should restrict the range of parameter $\theta \in (0, \pi)$ to the range $\theta \in [\theta_0, \theta_1]$ {where $\theta_0 = 2\cdot\arctan(1/e) \cong 0{,}2244\pi$, $e = 2.71828...$, $\theta_1 = 2\cdot\arctan(e) \cong 0{,}7756\pi$} for the reason that inequality: $|\ln \tan(\theta/2)| < |\tan(\theta/2)| \leq 1$ should be valid for all meanings of function $\ln(\tan(\theta/2))$ in that range of $\theta$, especially if $|k \cdot R| > \pi/4$.

So, these unusual beams *with limited* amplitude $A$ could be comparing to the spherical-waves solution (which is much more close to each other than other exotic beams) only at the range of parameter $\theta \in [{\sim}40{,}4°, {\sim}139{,}6°]$.

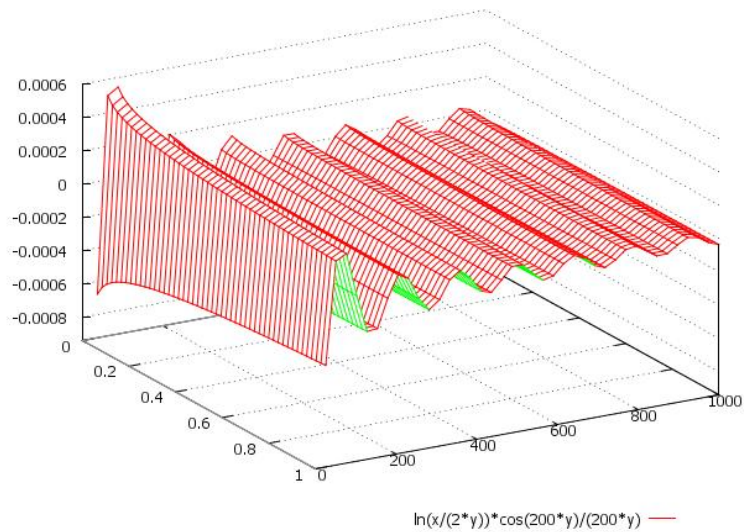

Fig.3. A *schematic* plot of the function $\ln(x/2y)*(\cos(k\cdot y)/(k\cdot y))$,
here we designate: $x = r(X,Y) = \sqrt{(X^2+Y^2)} \in (0, 1)$, $y = Z \in (0, 1000)$.

As for the appropriate examples of *paraxial* approximation $r(X,Y) = \sqrt{(X^2+Y^2)} \ll Z$, expressed by Eq. (4.1) in Cartesian co-ordinates *X, Y, Z*, see Fig.3,4.

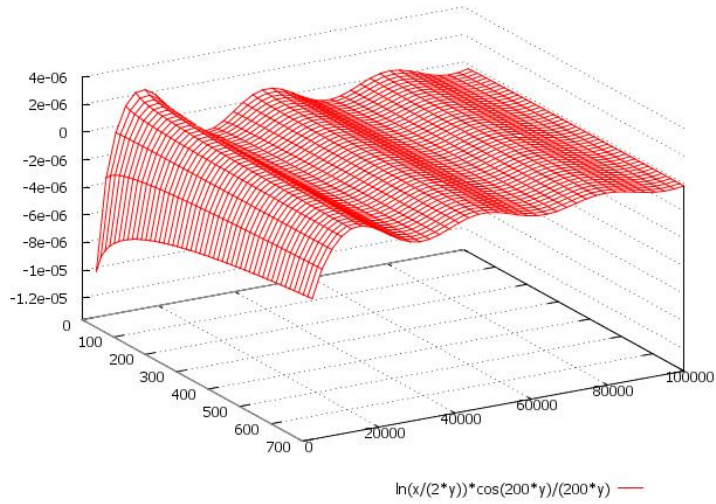

Fig.4. A *schematic* plot of the function $\ln(x/2y)*(\cos(k \cdot y)/(k \cdot y))$,

here we designate: $x = r(X,Y) = \sqrt{(X^2+Y^2)} \in (0, 700)$, $y = Z \in (0, 100'000)$.

Let us also *schematically* imagine the spherical-wave solution to compare it with the solution above:

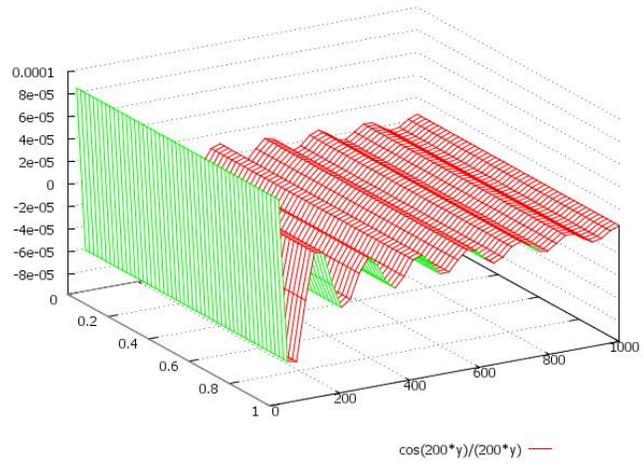

Fig.5. A *schematic* plot of the function $(\cos(k \cdot y)/(k \cdot y))$,

here we designate: $y = Z \in (0, 1'000)$.

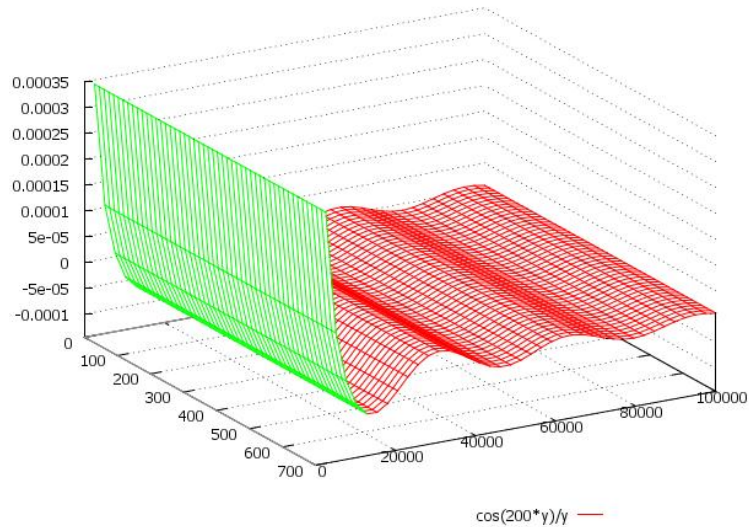

Fig.6. A *schematic* plot of the function (cos(k·y)/ k·y),

here we designate: $y = Z \in (0, 100'000)$.

Also, let us note that these unusual beams could be comparing to the Bessel beam solutions [11] at all the range of parameter $\theta \in (0, \pi)$. To obtain the energy of the beam, we should spherically integrate the absolute magnitude of the solution over space, so such a calculations should produce the infinite energy of a beam due to the structure of the solution: $A = (a \cdot k) \cdot \ln(\tan(\theta/2)) \cdot (\cos(k \cdot R))/(k \cdot R)$.

As for the point of clarifying the physical content of the derived solution, Fig. 3-6 could present a wave travelling on the ocean surface, for example.

Such a solution is supposed to be linearly enhancing their radius *R* during propagation in *R*-direction, but it is modulated by the function $\ln(\tan(\theta/2))$ so that the total amplitude $A = (a \cdot k) \cdot \ln(\tan(\theta/2)) \cdot (\cos(k \cdot R))/(k \cdot R)$.

Jumping of a phase-function of a solution (*) in a form (3.5) for an amplitude *A* being equal to zero at the meaning of parameter $\theta = \pi/2$, could be associated with the existence of an *optical vortex* [10] at this point. Optical vortex (also known as a screw dislocation or phase singularity) is a zero of an optical field, a point of zero intensity. Research into the properties of vortices has thrived since a comprehensive paper [13], described the basic properties of "dislocations in wave trains".

http://en.wikipedia.org/wiki/Gaussian_beam (see "Mathematical form").

11. P. W. Milonni and J. H. Eberly (2010). *Laser Physics* (2nd ed.). Section 14.14.

12. Carl M. Bender, Steven A. Orszag (1999). *Advanced Mathematical Methods for Scientists and Engineers*. Origin.published by McGraw Hill, 1978, XIV. pp. 20-22.

13. J. F. Nye, M. V. Berry (1974). *Dislocations in wave trains*. Proceedings of the Royal Society of London, Series A 336 (1605): 165. See also:

http://en.wikipedia.org/wiki/Optical_vortex

## Appendix (checking of the exact solution).

The direct substitution of the final expression (3.5) into the Helmholtz equation (1.1) is an easy matter, showing that this is really an exact solution.

Let us begin to check the solution (3.5) from the 1-st part of such a solution:

$$A = a \cdot \ln\left(\tan\left(\frac{\theta}{2}\right)\right) \cdot \frac{\cos(k \cdot R)}{R} = a \cdot \left(\int\left(\frac{1}{\sin\theta}\right) d\theta\right) \cdot \frac{\cos(k \cdot R)}{R} \quad \{|k \cdot R| < \pi/4, \ \theta \in (0, \pi)\}$$

For the reason some of readers may be have no the sufficient time to execute the calculations properly, it has been made step-by-step below:

$$\frac{\partial^2 A}{\partial R^2} + \frac{2}{R}\frac{\partial A}{\partial R} + k^2 A + \frac{1}{R^2}\left(\frac{\partial^2 A}{\partial \theta^2} + \cot\theta \frac{\partial A}{\partial \theta}\right) = 0, \quad \Rightarrow$$

$$a \cdot \left(\int\left(\frac{1}{\sin\theta}\right) d\theta\right) \cdot \left(\frac{d^2\left(\frac{\cos(k \cdot R)}{R}\right)}{dR^2} + \frac{2}{R}\frac{d\left(\frac{\cos(k \cdot R)}{R}\right)}{dR} + k^2\left(\frac{\cos(k \cdot R)}{R}\right)\right) +$$

$$+ \frac{1}{R^2} \cdot a \cdot \frac{\cos(k \cdot R)}{R} \cdot \left(\frac{d^2\left(\int\left(\frac{1}{\sin\theta}\right) d\theta\right)}{d\theta^2} + \cot\theta \frac{d\left(\int\left(\frac{1}{\sin\theta}\right) d\theta\right)}{d\theta}\right) = 0,$$

$$\Rightarrow \; a \cdot \left( \int \left( \frac{1}{\sin \theta} \right) d\theta \right) \cdot \left[ \frac{d^2 \left( \frac{\cos(k \cdot R)}{R} \right)}{dR^2} + \frac{2}{R} \frac{d \left( \frac{\cos(k \cdot R)}{R} \right)}{dR} + k^2 \left( \frac{\cos(k \cdot R)}{R} \right) \right] +$$

$$+ \frac{1}{R^2} \cdot a \cdot \frac{\cos(k \cdot R)}{R} \cdot \left( -\frac{\cos \theta}{(\sin \theta)^2} + \frac{\cos \theta}{\sin \theta} \left( \frac{1}{\sin \theta} \right) \right) = 0, \; \Rightarrow$$

$$\frac{d \left( \frac{-k \cdot \sin(k \cdot R) \cdot R - \cos(k \cdot R)}{R^2} \right)}{dR} + \frac{2}{R} \cdot \left( \frac{-k \cdot \sin(k \cdot R) \cdot R - \cos(k \cdot R)}{R^2} \right) + k^2 \left( \frac{\cos(k \cdot R)}{R} \right) = 0,$$

$$\frac{[-k^2 \cdot \cos(k \cdot R) \cdot R - k \cdot \sin(k \cdot R) + k \cdot \sin(k \cdot R)] \cdot R^2 - 2R \cdot [-k \cdot \sin(k \cdot R) \cdot R - \cos(k \cdot R)]}{R^4} +$$

$$+ \frac{2}{R} \cdot \left( \frac{-k \cdot \sin(k \cdot R) \cdot R - \cos(k \cdot R)}{R^2} \right) + k^2 \left( \frac{\cos(k \cdot R)}{R} \right) = 0, \; \Rightarrow$$

$$- \frac{k^2 \cdot \cos(k \cdot R)}{R} + \frac{2 \cdot k \cdot \sin(k \cdot R)}{R^2} + \frac{2 \cdot \cos(k \cdot R)}{R^3} +$$

$$+ \frac{k^2 \cdot \cos(k \cdot R)}{R} - \frac{2k \cdot \sin(k \cdot R)}{R^2} - \frac{2 \cos(k \cdot R)}{R^3} = 0,$$

- where the last identity is obviously valid for the range of parameter *R*: |*k·R*| < π/4.

The checking of 2-nd part of a solution (3.5) could be executed in the same way.